\documentclass[10pt,letterpaper]{article}
\usepackage[top=0.5in,left=1.5in,footskip=0.75in,marginparwidth=1in]{geometry}

\usepackage[utf8]{inputenc}

\usepackage{cite}

\usepackage{nameref,hyperref}

\usepackage{microtype}
\DisableLigatures[f]{encoding = *, family = * }

\raggedright
\usepackage{changepage}

\usepackage[aboveskip=1pt,labelfont=bf,labelsep=period,singlelinecheck=off]{caption}

\makeatletter
\renewcommand{\@biblabel}[1]{\quad#1.}
\makeatother

\usepackage{lastpage,fancyhdr,graphicx}
\usepackage{epstopdf}
\fancyheadoffset[L]{2.25in}
\fancyfootoffset[L]{2.25in}

\usepackage{color}

\definecolor{Gray}{gray}{.25}

\usepackage{graphicx}

\usepackage{sidecap}

\usepackage{wrapfig}
\usepackage[pscoord]{eso-pic}
\usepackage[fulladjust]{marginnote}
\reversemarginpar

\usepackage{bm}

\begin{document}
\vspace*{0.35in}

\begin{flushleft}
{\Large
\textbf\newline{Quadcubic interpolation: a four-dimensional spline method.}
}
\newline
\\
Walker, Paul\textsuperscript{1,*}
\\
\bigskip
\bf{1}  Department of Physics, Durham University, South Road, Durham, DH1 3LE, United Kingdom.

\bigskip
* paul.a.walker@durham.ac.uk

\end{flushleft}

\section*{Abstract}

We present a local interpolation method in four dimensions utilising cubic splines. An extension of the three-dimensional tricubic method, the interpolated function has C$^1$ continuity and its partial derivatives are analytically accessible. The specific example of application of this work to a time-varying three-dimensional magnetic field is given, but this method would work equally well for a time-independent four-dimensional field. Implementations of both of these methods in the Python programming language are also available to download.




\section*{Introduction}

Where some quantity of interest, for example magnetic field strength or temperature, is known as a discrete data set, it is necessary to interpolate the set to acquire knowledge of the value at arbitrary points within the data domain. Various schemes exist, such as `polynomial interpolation' which fits some high-level function to all data points. These can be computationally expensive to calculate, and can suffer from oscillations, particularly at the boundaries, behaviour known as Runge's phenomenon \cite{runge}. Alternatively there are piecewise methods, in which some approximating function is found between adjacent data points. A simple example of this is linear interpolation, but this does not produce a very smooth output - better are splines \cite{ahlberg}, which use some low-degree (typically cubic) polynomial function. These do not suffer from Runge's phenomenon, are relatively cheap to calculate, are smoother than linear methods and, being local, large errors (often found at data boundaries) do not propagate.\\
For systems with more than one variable multivariate schemes exist to combine interpolation of the data set across these variables - for example, two-dimensional linear and cubic methods are readily available in Scientific Python \cite{scipyinterp}. Originally motivated by studies of ocean dynamics \cite{tricubic2} Lekien and Marsden describe their `tricubic' technique \cite{tricubic} for either time-dependent two-dimensional flows or three-dimensional time-independent flows. Earlier three-dimensional cubic spline methods treated the problem as three one-dimensional problems \cite{catmull}, whereas the Lekien-Marsden solution efficiently combines them.\\
We implemented this interpolator in Python \cite{walkertricubic} for use in modelling the motion of paramagnetic neutral particles through Zeeman decelerators \cite{mcard} and magnetic traps \cite{walker}. These fields are produced by combinations of permanent magnetic and electromagnetic elements with generally no analytic solution. The potentials are calculated, for example using finite element analysis methods, as a series of data points on a grid, which must be interpolated to return the required values for an arbitrary point within the region of interest. In order to more effectively incorporate time-varying magnetic fields we have extended this model to be four-dimensional.\\

\section*{Tricubic interpolation}
The tricubic method given by Lekien and Marsden \cite{tricubic} assumes a data set known at a series of regular grid points in three dimensions $x,y$ and $z$. As long as the data are regular in all three dimensions they do not have to be the same length, but without loss of generality we shall here consider a coordinate mesh that is a cube. This mesh is composed of elements that are also cubes, for each of which the interpolant field $f$ is known at the eight vertices $p_1 \cdots p_8$ (this is, of course, for a scalar field - for a vector field the same method is simply applied to each component separately). Inside each element $f$ is in the form of the cubic:

\begin{equation}\label{eq:tricubic}
f(x,y,z) = \sum_{i,j,k=0}^{3} a_{ijk} x^i y^j z^k.
\end{equation}

\noindent The coefficients $a_{ijk}$ must be calculated, and can then be used to return $f$ for an arbitrary point $(x,y,z)$ within the element. In order to achieve $C^1$ continuity across the whole domain, the values of $f$ and its first derivatives $\partial f / \partial x, \partial f / \partial y$ and $\partial f / \partial z$ must be continuous across the faces and vertices of the elements, giving 32 constraints for each. The list of coefficients in equation \ref{eq:tricubic} is of length 64, so there must be an additional 32 constraints; it can be shown \cite{tricubic} that the only valid choice is the set:

\begin{equation}\label{eq:triconstraints}
\left[ f, \frac{\partial f}{\partial x}, \frac{\partial f}{\partial y}, \frac{\partial f}{\partial z}, \frac{\partial^2 f}{\partial x\partial y}, \frac{\partial^2 f}{\partial x\partial z}, \frac{\partial^2 f}{\partial y\partial z}, \frac{\partial^3 f}{\partial x\partial y \partial z} \right]
\end{equation}

\noindent To see why derivatives such as $\partial^2 f / \partial x^2$ cannot be used, consider the value of $\partial^2 f / \partial x^2$ at point $p_1$, which we can take to be at the origin $(0,0,0)$. Along the $x-$axis equation \ref*{eq:tricubic} reduces to the cubic spline:

\begin{equation}\label{eq:xcubic}
f(x) = \sum_{i=0}^{3} a_{i} x^i.
\end{equation}

\noindent The values of $f$ and $\partial f / \partial x$ at points $p_1$ and $p_2 = (1,0,0)$ are fixed, forming a unique spline. This also constrains $\partial^2 / \partial x^2$ at $p_1$ and $p_2$, which is therefore not independent of the values of $f$ and $\partial f / \partial x$. More generally this applies to all derivatives of the form $\partial^2 / \partial z^2$, $\partial^3 / \partial y \partial z^2$, \textit{etc.}, and means that no cubic spline method can achieve C$^2$ continuity \cite{press}.\\

The derivatives in set \ref{eq:triconstraints} are found via finite-difference methods \cite{zhou}. If the coefficients $a_{ijk}$ from equation \ref{eq:tricubic} are placed into a vector $\bm{\alpha}$ and the values of $f$ and its derivatives from set \ref{eq:triconstraints} into another vector $\bm{b}$, they are related by a 64 $\times$ 64 matrix $B$:

\begin{equation}
B \bm{\alpha} = \bm{b},	
\end{equation}

\noindent where the components of $B$ are integers. If the elements are normalised to be unit cubes during calculations, this matrix is the same for all elements, and so only needs to be calculated once. During interpolation the results are scaled back to their actual values. This matrix is too large to include here but we have made it available online \cite{walkertricubic}. The matrix is invertible, so:

\begin{equation}
	B^{-1} \bm{b} = \bm{\alpha},
\end{equation}

\noindent allowing the coefficients $\bm{\alpha}$ to be calculated. In order to perform an interpolation for an arbitrary point $(x,y,z)$, firstly the appropriate volume element containing the query point is located. The alpha coefficients for this element are calculated, and then combined with the query coordinates in equation \ref{eq:tricubic} to return $f$. The coefficients $a_{ijk}$ are a tensor of order 3, which is applied to the vectors $(1, x, x^2, x^3), (1, y, y^2, y^3)$ and $(1, z, z^2, z^3)$. Differentiating these vectors can be used to return the derivates of $f$ - for example, using the set $(0, 1, 2x, 3x^2), (1, y, y^2, y^3)$ and $(1, z, z^2, z^3)$ gives $\partial f / \partial x$. We recently released \cite{walkertricubic} an implementation of this tricubic interpolator written in Python \cite{python} and NumPy \cite{numpy}. The online repository also contains example scripts and data files.

\section*{Quadcubic interpolation}
Motivated by a requirement to model particle trajectories through time-dependent magnetic fields, we have extended the Lekien Marsden method to four dimensions. Instead of a cube in $(x,y,z)$ with eight vertices, each volume element in our interpolation region is a tesseract in $(x,y,z,t)$ with 16 vertices, inside of which we have the function:

\begin{equation}\label{eq:quadcubic}
f(x,y,z,t) = \sum_{i,j,k,l=0}^{3} a_{ijkl} x^i y^j z^k t^l.
\end{equation}

\noindent There are 256 coefficients $a_{ijkl}$ meaning 256 constraints are required, and the derivatives given in set \ref*{eq:triconstraints} only supply 128 of them. Again noting the requirement to add independent constraints we add the additional derivatives in set \ref*{eq:quadconstraints} to set \ref*{eq:triconstraints}:

\begin{equation}\label{eq:quadconstraints}
	\left[\frac{\partial f}{\partial t}, \frac{\partial^2 f}{\partial x\partial t}, \frac{\partial^2 f}{\partial y\partial t}, \frac{\partial^2 f}{\partial z\partial t}, \frac{\partial^3 f}{\partial x\partial y \partial t}, \frac{\partial^3 f}{\partial x\partial z \partial t}, \frac{\partial^3 f}{\partial y\partial z \partial t}, \frac{\partial^4 f}{\partial x\partial y \partial z\partial t} \right]
\end{equation}

\noindent In order to calculate the derivatives using finite-differences, the values of $f$ at the vertices of neighbouring elements is needed - a point in this space is surrounded by $3^4 = 81$ tesseracts defined by 256 points (as many of the vertices are shared). Placing these in a vector $\bm{x} = c_1 \cdots c_{256}$, the finite-difference matrix \cite{zhou} $D$ returns the elements of set \ref{eq:triconstraints} as a vector $\bm{b}$:

\begin{equation}
D \bm{x} = \bm{b}.
\end{equation}

\noindent The interpolation matrix $B$ is 256 $\times$ 256 in size. Additional savings in computation time can be achieved by following the method of Faust \textit{et al.} \cite{eqtools}; we can combine the matrices $D$ and $B^{-1}$ to produce a new matrix $D B^{-1} = A$, which only needs to be calculated once and then reused as required:

\begin{equation}
A \bm{x} = \bm{\alpha}
\end{equation}

\noindent $a_{ijkl}$ is a tensor of order 4, and we apply it to the four-dimensional vectors  $(1, x, x^2, x^3), (1, y, y^2, y^3)$, $(1, z, z^2, z^3)$ and $(1, t, t^2, t^3)$ to return $f$. As before, taking the derivatives of one of these vectors with respect to its variable allows the interpolated values of $\partial f / \partial x$, $\partial f / \partial t$ \textit{etc.} to be calculated.\\
This new interpolation method has been added to the `ARBInterp' library \cite{walkertricubic} and is available to download and use on a GPL open-source license. Once again, thanks to Dr. Lewis McArd for his invaluable advice on this and other projects. This work was undertaken as part of research funded by EPSRC grant number EP/N509462/1.

\end{document}